\documentclass{svproc}
\usepackage{url}

\usepackage{amsmath,amssymb}
\usepackage{tikz-cd}

\DeclareMathOperator{\spn}{span}

\begin{document}
\mainmatter              
\title{A tutorial on the dynamic Laplacian}
\titlerunning{A tutorial on the dynamic Laplacian}  
%
\author{Gary Froyland}
\authorrunning{Gary Froyland} 
%
\tocauthor{Gary Froyland}
\institute{UNSW Sydney, NSW 2052, Australia\\
\email{g.froyland@unsw.edu.au}\\
\texttt{www.maths.unsw.edu.au/$\sim$froyland}\\
\texttt{www.dynamicdata.unsw.edu.au}}

\maketitle              

\begin{abstract}
Spectral techniques are popular and robust approaches to data analysis. 
A prominent example is the use of eigenvectors of a Laplacian, constructed from data affinities, to identify natural data groupings or clusters, or to produce a simplified representation of data lying on a manifold.
This tutorial concerns the dynamic Laplacian, which is a natural generalisation of the Laplacian to handle data that has a time component and lies on a time-evolving manifold.
In this dynamic setting, clusters correspond to long-lived ``coherent'' collections.
We begin with a gentle recap of spectral geometry before describing the dynamic generalisations.
We also discuss computational methods and the automatic separation of many distinct features through the SEBA algorithm.
The purpose of this tutorial is to bring together many results from the dynamic Laplacian literature into a single short document, written in an accessible style.

%
\end{abstract}
\section{Introduction}

Spectral geometry is concerned with relationships between the geometry of manifolds on the one hand, and spectral information from differential operators on manifolds on the other, most commonly the Laplace operator.
While so-called direct problems seek spectral information from manifold geometry, inverse problems attempt to learn manifold geometry from the operator spectrum.
Inverse problems are of great interest in practical areas such as spectral approaches to dimension reduction of high-dimensional manifolds and geometric clustering of data.

This brief tutorial focusses on the inverse problem, and its generalisation to answer questions concerning long-lived structurs in dynamical systems.
In this dynamic setting, the manifold is subjected to nonlinear transformations representing the evolution of the manifold over time.
The inverse problem in this setting concerns learning the average geometry of the evolved manifold.
Corresponding to this dynamic inverse problem is a natural dynamic generalisation of the Laplace operator to a \emph{dynamic Laplace operator}.
Using both the spectrum and eigenfunctions of the dynamic Laplacian we create a dynamic spectral geometry, and illustrate its use in identifying \emph{finite-time coherent sets} in dynamical systems.
These sets are the natural generalisations of Cheeger sets to the dynamic setting.


\section{Recap on properties of the Laplace operator}
\label{sec:laplace}
We begin by recapping some basic spectral geometry for manifolds in the absence of dynamics;  see for example \cite{lablee}.
Consider a full-dimensional manifold\footnote{All of the properties discussed here hold on general Riemannian manifolds, but for simplicity of presentation we will consider flat manifolds $M$ endowed with the Euclidean metric.} $M\subset\mathbb{R}^d$.
Expressed in local coordinates, the Laplace operator $\Delta$ acts on functions $f:M\to\mathbb{R}$ by $\Delta f=\sum_{i=1}^d \frac{\partial^2}{\partial x_i^2}f$.
The eigenproblem is
\begin{equation}
\label{lapeig}
\Delta f = \lambda f,
\end{equation}
with $\langle \nabla f(x), n(x)\rangle=0$ for all $x\in\partial M$, where $n(x)$ is a vector normal to the boundary of $M$ at $x\in\partial M$, and $\langle\cdot,\cdot\rangle$ is the usual Euclidean inner product.
With these Neumann boundary conditions, the eigenvalues 
 consist of a countable sequence $0=\lambda_1\ge \lambda_2\ge \cdots$, of negative numbers, diverging to $-\infty$.
The leading eigenvalue $\lambda_1=0$ and the corresponding eigenfunction is constant on $M$.
The eigenfunctions $f_1,f_2,\ldots$ form a complete orthonormal set for $L^2(M)$.
\subsection{Optimality properties of Laplace eigenvalues and eigenfunctions}
\begin{figure}
\centering
  \includegraphics[width=\textwidth]{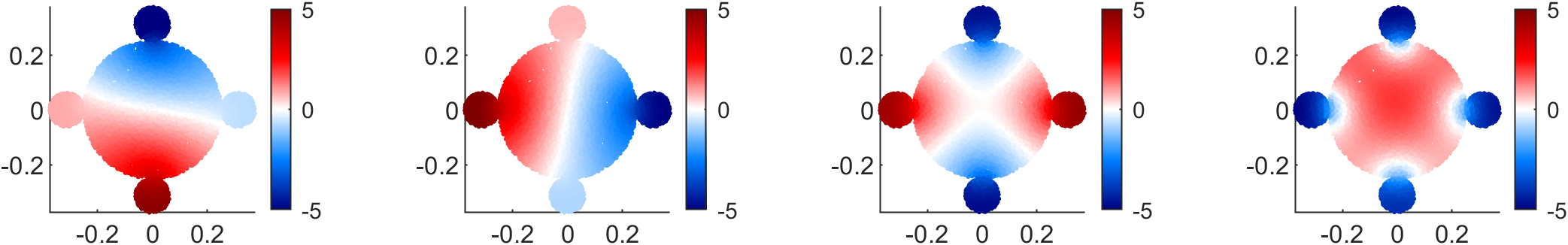}
  \caption{Laplace eigenfunctions $f_2,\ldots,f_5$ of a two-dimensional manifold. The corresponding eigenvalues are $\lambda_2\approx -33.42$, $\lambda_3\approx -33.80$, $\lambda_4\approx -49.99$, $\lambda_5=-82.42$.}\label{fig:blob}
\end{figure}

There are variational expressions for the Laplace eigenvalues, for example, for each $n\ge 1$,
\begin{equation}
\label{variational}
\sum_{k=1}^n \lambda_k = \min\left\{\sum_{k=1}^n \int_M |\nabla f_k(x)|^2\ dx : f_1,\ldots,f_n\mbox{ are orthonormal}\right\}.
\end{equation}
Moreover, for each $n\ge 1$, the minimizing collection $\{f_1,\ldots,f_n\}$ is precisely the $n$ leading Laplace eigenfunctions.
Thus, denoting $L_n=\spn\{f_1,\ldots,f_n\}$, we see that the Laplace eigenfunctions are the most regular (with respect to the $L^2$ norms of gradients) orthonormal basis for $L_n$.
Equivalently to \eqref{variational} we have for each $n\ge 1$
\begin{equation}
\label{variational2}
 \lambda_n = \min\left\{\int_M |\nabla f_n(x)|^2\ dx : \|f_n\|_2=1, f_n\perp L_{n-1}\right\},
\end{equation}
where $L_0$ is the trivial subspace.
Figure \ref{fig:blob} shows the first four nontrivial eigenfunctions of the Laplace operator on a two-dimensional manifold consisting a large disk with four smaller disks connected near the boundary of the large disk.
The decreasing regularity -- left to right in Figure \ref{fig:blob} -- as indicated in (\ref{variational}) and \eqref{variational2} is apparent.

\subsection{The heat equation and geometry}
The reader may also have noticed relationships between the eigenfunctions in Figure \ref{fig:blob} and the geometry of the manifold.
Indeed there are well-known links between the Laplace eigendata of a manifold and the geometric shape of the manifold.
Firstly, consider the heat equation on a manifold $M$, where the scalar value $u(t,x)$ denotes the temperature at position $x$ and time $t$:
\begin{equation}
\label{heateqn}
\frac{\partial u}{\partial t}=\Delta u.
\end{equation}
We impose the same Neumann boundary conditions as above;  this means that there is no heat lost at the boundary due to a zero derivative in the temperature field normal to the boundary.
If we initialise a heat flow with $u(0,x)=f_k(x)$ then we obtain $u(t,x)=e^{\lambda_k t}f_k(x)$, $k=1,2,\ldots$.
The minimality properties \eqref{variational} of the eigenvalues $\lambda_k$ corresponds to a minimal decay of an orthonormal set of initial conditions $f_1,\ldots,f_n$ toward the heat equilibrium of $f_1\equiv 1/\sqrt{\mbox{vol}(M)}$. 
In particular, in the two left panels of Figure \ref{fig:blob} we see that positive ``heat'' is concentrated to one side of the manifold and negative heat is concentrated on the other side, with the greatest concentration in the opposing small disks. Such an initialisation will lead to a lower heat relaxation toward equilibrium.
Similarly, in the right two panels of Figure \ref{fig:blob} the heat is initialised predominantly in the small disks around the periphery of $M$.
Because of the narrow necks joining the small disks to the main disk the heat flow through these necks will contribute to slow relaxation to equilibrium.

\subsection{Cheeger problems and Laplace eigenfunctions}
\label{sec:cheelap}
Similar notions can be constructed in a purely geometric fashion through \emph{Cheeger problems} \cite{cheeger,buser}.
Let $\mathring M$ denote the interior of $M$ and $|A|$ denote the Lebesgue volume of $A$, matched to the dimension of the set. 
For a set $A\subset M$ with smooth boundary, we call 
\begin{equation}
\label{neumcheegerratio}
H(A):=\frac{|\partial A \cap \mathring M|}{|A|}
\end{equation} the \emph{Neumann Cheeger ratio} of $A$, which quantifies the boundary size of $A$ relative to its interior size.
We call a collection $A_1,\ldots,A_n$ of $n$ disjoint subsets of $M$ an \emph{$n$-packing} and define $H(\{A_1,\ldots,A_n\})=\max_{1\le k\le n}H(A_k)$.
Cheeger (or isoperimetric) problems are concerned with decomposing manifolds into subsets with small boundaries relative to their interiors.
To this end, we define the \emph{$n^{\rm th}$ Neumann Cheeger constant} of $M$ as 
\begin{equation}
\label{cheeger}
h_{n}:=\inf \{H(\{A_1,\ldots,A_n\}): A_1,\ldots,A_n\subset M\mbox{ are disjoint}\};
\end{equation}
see \cite{FR24}.
It is not hard to show that $h_n$ is increasing in $n$.
Cheeger inequalities for $n$-packings of boundaryless manifolds have been studied in e.g.\ \cite{miclo15}, with analogous definitions.

In general, it is difficult to find the minimising collection $A_1,\ldots,A_n$, however the Laplace eigenfunctions can provide good candidate collections.
A \emph{nodal domain} of $f_k$ is a maximal connected component of $M$ where the sign of $f_k$ is constant.
We have the following result
\begin{theorem}[Theorem 3.7 \cite{FR24}]
\label{ncheegerthm}
Consider an eigenpair $(\lambda,f)$.
\begin{enumerate}
\item If $f$ has at least $n$ nodal domains then $h_{n}\le \sqrt{-2\lambda}$.
\item For each $k=1,\ldots,n$, let $N_k\subset M$ denote the $k^{\rm th}$ nodal domain and let $A_k(c_k):=\{x\in N_k: f(x)^2\ge c_k\}\subset M$ denote the superlevel set of $f^2$ restricted to $N_k$.
    Then for a positive measure set of $(c_1,\ldots,c_n)\in\mathbb{R}^n$, one has 
    \begin{equation}
    \label{feassol}
    H(A_1(c_1),\ldots,A_n(c_n))\le \sqrt{-2\lambda}.
    \end{equation}
\end{enumerate}
\end{theorem}
Part 1 says that the closer $\lambda$ is to zero, the more easily $M$ can be decomposed into $n$ pieces, each with a small Cheeger ratio.
This is consistent with our discussion of the heat equation because the closer $\lambda$ is to zero, the more slowly the initial condition $u(0,x)=f(x)$ will relax to equilibrium. 
When $f$ has $n$ nodal domains, there are $n$ subsets of $M$  containing obstacles that slow down this equilibriation.
In Figure \ref{fig:blob}, the left two panels show that $f_2$ and $f_3$ have two nodal domains (separated by the white strip), each roughly half of $M$, $f_4$ has four nodal domains, and $f_5$ has five nodal domains.
Thus, by part 1 of Theorem \ref{ncheegerthm}, using the eigenvalues quoted in the caption of Figure \ref{fig:blob} we have $h_2\le 8.18$, $h_3\le h_4\le 10$, and $h_5\le 12.84$.

Part 2 of Theorem \ref{ncheegerthm} provides a recipe for constructing $n$-packings of $M$ with low Cheeger ratios, namely level sets of $f^2$ individually restricted to each nodal domain.
In this regard, note that the nodal lines of $f_2$ and $f_3$ in Figure \ref{fig:blob} are slightly offset from the horizontal and vertical, respectively;  this is not an accident.
If the nodal lines were horizontal and vertical, they would split the small disks at the periphery of $M$,  leading to larger Cheeger ratios for the corresponding blue and red nodal domains.
Cheeger inequalities without nodal domain counts are available for boundaryless manifolds \cite{miclo15}, but they are much weaker.

\subsection{Sparse eigenbasis approximation (SEBA)} 
\label{sec:seba1}
Finding good candidate $n$-packings of a manifold with superlevel sets of individual nodal domains of $f^2$ can be practically challenging.
Computing \emph{individual} nodal domains can be nontrivial because of the connectivity requirement in the definition.
Moreover, because the irregularity of the $\{f_k\}$ grows with $k$, we expect superior $n$-packings to be obtained from as few $f_k$ as possible.
The $n^{\rm th}$ eigenfunction $f_n$ contains \emph{at most} $n$ nodal domains (by the Courant Nodal Domain Theorem \cite{courant}), and may contain fewer.

\begin{figure}
\centering
  \includegraphics[width=\textwidth]{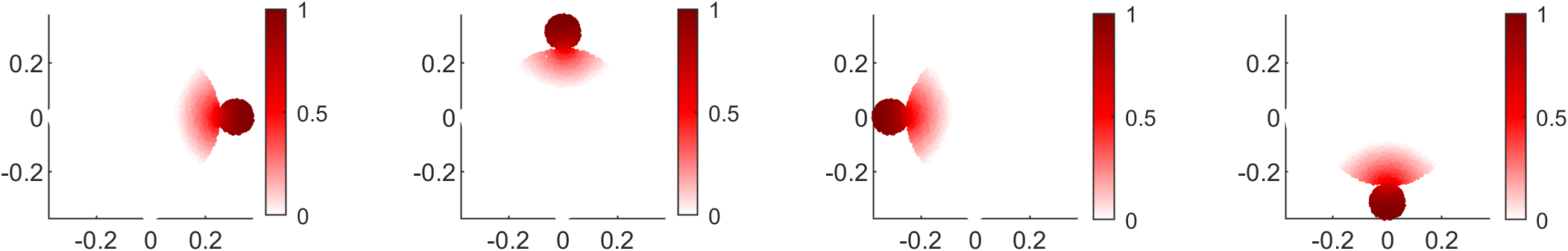}
  \caption{Sparse functions $s_1,\ldots,s_4$ output by the SEBA algorithm applied to $f_1,\ldots,f_4$. One has $\spn\{s_1,\ldots,s_4\}\approx \spn\{f_1,\ldots,f_4\}$ and may define $A_k$ as a superlevel set of $s_k$, $k=1,\ldots,4$.}\label{fig:blobSEBA}
\end{figure}

By taking linear combinations of $\{f_1,\ldots,f_n\}$ in combination with soft thresholding we can create a sparse collection of functions $\{s_1,\ldots,s_n\}$ that approximately span $\spn\{f_1,\ldots,f_n\}$. 
Each function $s_k$ will support an individual element of an $n$ packing;  see \cite{FR24} for a discussion of SEBA in the manifold setting, including a direct analogue to Part 2 of Theorem \ref{ncheegerthm} for sets $A_k(c_k):=\{x\in M: s_k(x)\ge c_k\}$.
Implementational details, further discussion of constructing $n$-packings, and the short code may be found in \cite{FRS19}.
Figure \ref{fig:blobSEBA} illustrates this approach using eigenfunctions $f_1,\ldots,f_4$ (i.e.\ the constant function $f_1$ along with the left 3 panels of Figure \ref{fig:blob}).
SEBA has worked well when the sparse basis functions $s_1,\ldots,s_n$ are all nonnegative. 
For example, applying SEBA to the collections of eigenfunctions $\{f_1,\ldots,f_n\}$, $2\le n\le 5$ in Figure \ref{fig:blob}, we obtain nonnegative basis functions for $n=2,4,5$, but not $n=3$ because there is no natural encoding of three sets in the leading three eigenfunctions.

Finally, we remark that the $k$-means algorithm is often used to post-process spectral clustering output.
Our Laplacian eigenanalysis is precisely spectral clustering on the manifold $M$ and SEBA provides a more robust post-processing step for discrete data than $k$-means.
This is because SEBA will zero out data points that do not naturally belong to a cluster and SEBA provides clear guidelines for rejecting spurious sparse vectors ($s_k$ with large negative values should be rejected);  see \cite{FRS19} for details.


\section{Dynamical systems and finite-time coherent sets}
We now introduce dynamics on manifolds and explain how dynamic versions of the spectral geometric properties of the Laplacian discussed in the previous section can identify coherent sets in dynamical systems.
We consider a domain\footnote{In dynamical systems, the domain is often called phase space or state space.} $M$ and one-parameter family\footnote{All of the constructions developed here have obvious analogues if time is discrete.} of smooth transformations $\Phi^t:M\to\Phi^t(M)$ denoting the evolution of points in $M$ for a duration $t\in[0,\tau]$.
For example, $\Phi^t$ could be the flow generated by a time-dependent ODE $\dot{x}=F(t,x)$,\ i.e.\ $\frac{d}{ds}\Phi^s|_{s=t}(\cdot)=F(t,\cdot).$
For simplicity of presentation we will continue to assume that $M$ is $d$-dimensional compact submanifold of $\mathbb{R}^d$ endowed with the Euclidean metric, and that $\Phi^t$ preserves Euclidean volume.
All of the constructions below naturally generalise to evolving families of Riemannian manifolds endowed with metrics independent of the dynamics, and to non-Riemannian-volume-preserving dynamics $\Phi^t$ \cite{FK20}.

\subsection{Chaotic advection, mixing, and coherence}
\label{sec:chaoticadvection}
\begin{figure}
  \centering
  \includegraphics[width=\textwidth]{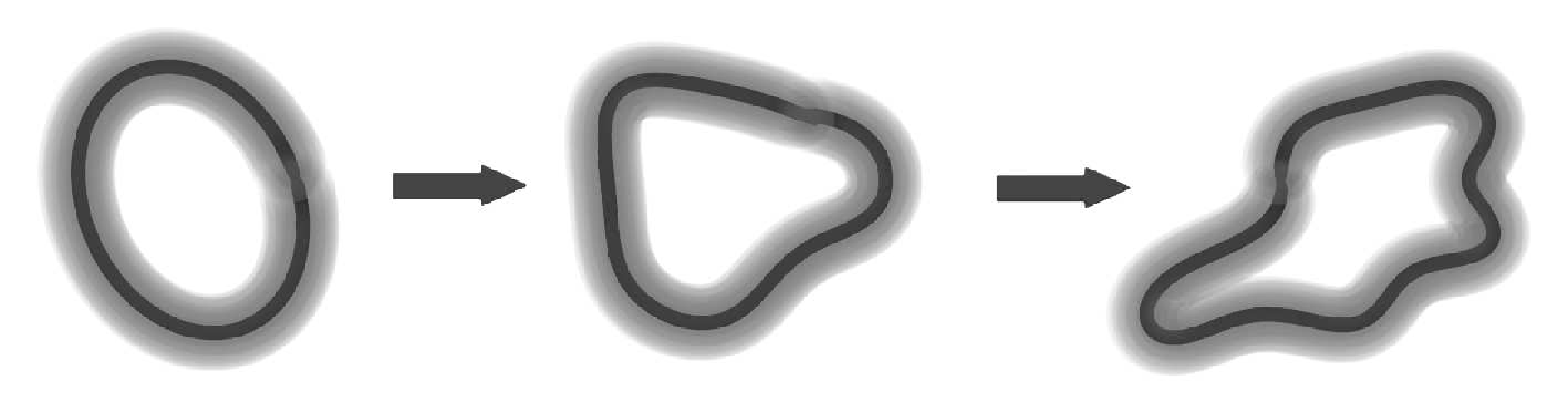}
  \caption{Advective mixing in the presence of small diffusion. Left: a two-dimensional boundary shown as a black curve with $\epsilon$-scale band of diffusion over the boundary. Centre and Right: under advective evolution with $\Phi^t$ the boundary evolves in a deterministic way, with the amount of mixing proportional to the size of the boundary.}\label{fig:diffusion}
\end{figure}

The above finite-time evolution is often highly heterogeneous in space, and an important characterisation of the dynamics is the decomposition of the phase space into coherent regions that are relatively non-dispersive and mixing regions that rapidly disperse.
We will say that a subset $A$ of $M$ is \emph{coherent over the interval $[0,\tau]$} if the boundary of $A$ remains small relative to its volume.
In the presence of small diffusion -- for example considering an SDE with drift $F$ as above and small diffusion,  
$dX_t=F(t,X_t)+(\epsilon^2/2)dW_t$ -- the points in the interior of $A$ that escape from $A$ are those that are near to the boundary, and any escape is by diffusion alone; see Figure \ref{fig:diffusion}.
This approach is consistent with methods of quantification of advective mixing \cite{ottino,pierrehumbert,petzold,thiffeault,haller2013}.

\begin{figure}
  \centering
  \includegraphics[width=\textwidth]{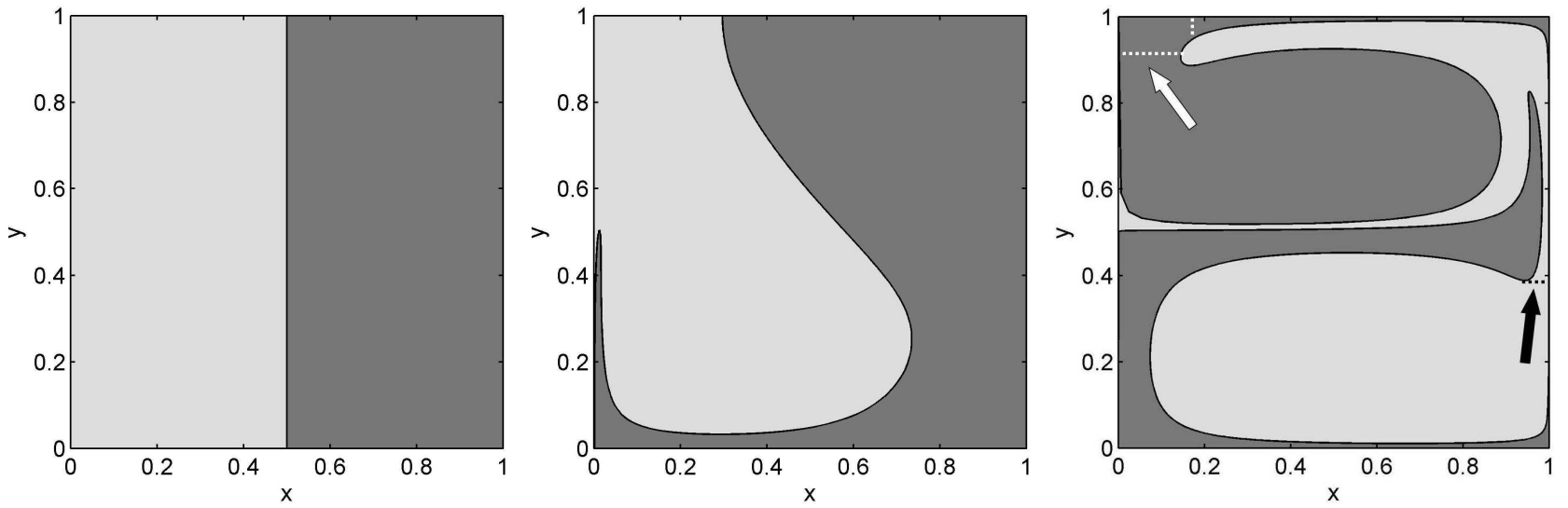}
  \caption{Illustration of the dynamics of the rotating double gyre flow. Left: at time $t=0$ we initialise with two sets $A_1, A_2$ given by the left and right halves of the domain. Centre: at time $t=1/2$ the boundaries of $A_1, A_2$ begins to grow. Right: at the final evolution time $t=1$ the boundaries have grown significantly, leading to advective mixing in the neighbourhoods of the boundaries. The arrows indicate asymmetries in the evolved sets. Image from \cite{FPG14}.}\label{fig:meiss_triple}
\end{figure}

We briefly illustrate this idea in a simple example, based on the rotating double gyre \cite{mosovskymeiss}. 
Let $M=[0,1]^2$ and $\Phi^t$, $t\in[0,1]$ be generated by the time-dependent ODE given by equation (3.1) \cite{mosovskymeiss}.
Figure \ref{fig:meiss_triple}(left) shows two sets $A_1,A_2\subset M$ with small Neumann Cheeger ratio (see \eqref{neumcheegerratio}): $H(A_1)=H(A_2)=1/(1/2)=2$. 
In fact, these $A_1, A_2$ form an optimal $2$-packing of $M$, with $h_2=2$.
But when these sets are evolved under the dynamics $\Phi^t$, their boundaries grow significantly: $H(\Phi^1(A_1))=H(\Phi^1(A_2))\approx 8.3057/(1/2)=16.6114$.
Thus, in a short period of time, exchange of fluid between the two shaded regions in Figure \ref{fig:meiss_triple} rapidly increases.

\subsection{Dynamic Cheeger ratios and coherent sets}
In view of Figure \ref{fig:diffusion} and the above discussion on advective mixing, coherent subsets $A$ of $M$ are those for which the total exchange of trajectories from $A$ to its complement, relative to the volume of $A$, is minimal.
This total exchange is proportional to the average Cheeger ratio of $\Phi^t(A)$ over the time interval $[0,\tau]$.
We therefore introduce a \emph{dynamic Cheeger ratio} \cite{F15,FR24} in analogue to the Neumann Cheeger ratio \eqref{neumcheegerratio}:
\begin{equation}
\label{dyncheegerratio}
H^D(A):=\frac{1}{\tau}\int_0^\tau \frac{|\partial(\Phi^t(A))\cap \Phi^t(\mathring{M})|}{|\Phi^t(A)|}\ dt.
\end{equation}
For $A=A_1$ in Figure \ref{fig:meiss_triple}, the numerator of \eqref{dyncheegerratio} is simply the average length of the common boundary of $A_1$ and $A_2$;  in particular, we do not include the boundaries\footnote{In Section \ref{sec:dirichlet} we will discuss the Dirichlet dynamic Cheeger ratio and the circumstances in which to use each of the two types of dynamic Cheeger ratios.} of $\Phi^t(M)$.
We therefore seek finite-time coherent sets as minimisers of $H^D(\cdot)$.
Our approach will mirror our search for minimisers of $H(\cdot)$ in Section \ref{sec:laplace}.
Directly minimising the geometric expression \eqref{dyncheegerratio} is difficult, so we instead develop a theory of \emph{dynamic spectral geometry}, and find good candidates for minimisers with eigenfunctions of a \emph{dynamic Laplace operator}.

\subsection{Encoding evolved geometries with pullback metrics}
\label{sec:pullback}
If our dynamic Laplace operator is going to help us find finite-time coherent sets, it must incorporate the deformed geometry of $M$ under the action of $\Phi^t$.
For example, in the right panel of Figure \ref{fig:meiss_triple} we computed the length of the evolved common boundary (let us denote this family of curves $\Phi^t(\Gamma$)) -- with respect to the Euclidean metric -- to be about $8.3057$ at time $t=1$.
To encode this deformation using a single operator acting on functions on $M$, we need to create a geometry on $M$ that reflects the Euclidean distances between points in $\Phi^t(M)$.
Under such a geometry, the length of $\Phi^t(\Gamma)$ -- computed with the Euclidean metric $e$ on $\Phi^t(M)$ -- would be the length of $\Gamma$ on $M$.
This geometry is described by the pullback metric $(\Phi^t)^*e=:g_t$.
Because we have assumed that $\Phi^t$ preserves volume on $M$, the metric $g_t$ has a particularly simple pointwise description in standard coordinates, namely it can be represented by the $d\times d$ positive-definite matrix $\mathbf{g}_t(x):=(D\Phi^t(x))^\top D\Phi^t(x)$ for $x\in M$, which is known as the Cauchy--Green tensor in fluid dynamics.


Continuing with the example in Figure \ref{fig:meiss_triple}, let $\gamma(s)=(1/2,s)$, $0\le s\le 1$ be a parameterisation of $\Gamma$ in the square domain $M$.
The length of $\Gamma$ according to the Euclidean metric on the square is 
\begin{equation}
\label{lineeqn}
\int_0^1 \|\gamma'(s)\|_2\ ds=\int_0^1 1\ ds=1.
\end{equation}
We may similarly compute the length of $\Phi^t(\Gamma)$.
A little notation before beginning: 
on the tangent space of $M$ at $x\in M$, denoted $T_xM$, we define an inner product $\langle \cdot,\cdot\rangle_{g_t}$ by $\langle u,v\rangle_{g_t}=u^\top\mathbf{g}_t(x)v$ for $u,v\in T_xM$.
We now may calculate the length of $\Phi^t(\Gamma)$ as:
\begin{eqnarray}
\nonumber\int_0^1 \|D\Phi^t(\gamma(s))\cdot \gamma'(s)\|_2\ ds&=&\int_0^1 \left(\gamma'(s)^\top(D\Phi^t(\gamma(s))^\top D\Phi^t(\gamma(s))) \gamma'(s)\right)^{1/2}\ ds\\
\nonumber&=&\int_0^1\left(\gamma'(s)^\top \mathbf{g}_t(\gamma(s)) \gamma'(s)\right)^{1/2}\ ds\\
\nonumber&=&\int_0^1 \langle \gamma'(s),\gamma'(s)\rangle_{g_t}^{1/2}\ ds\\
&=&\int_0^1 \|\gamma'(s)\|_{g_t}\ ds.
\label{lineeqn2}
\end{eqnarray}
Note that \eqref{lineeqn2} is almost identical in form to the left-hand-side of \eqref{lineeqn};  we have merely \emph{changed the norm} in which we compute the length of the derivative vector $\gamma'$.

As discussed above, the metric $g_t$ is the \emph{pullback metric} we obtain by pulling back the Euclidean metric $e$ with the diffeomorphism $\Phi^t$.
As Riemannian manifolds, $\Phi^t$ is an isometry between $(M,g_t)$ and $(\Phi^t(M),e)$, because lengths computed in the tangent space $T_xM$ with $g_t$ are identical to lengths computed in the tangent space $T_{\Phi^tx}\Phi^t(M)$ with the Eulidean metric.
The concept of pullback metrics allows us to \emph{transfer geometry on $\Phi^t(M)$ back onto $M$ by applying a new metric on $M$.}

If $M$ is $d$-dimensional, $d\ge 2$, we need to compute areas of $d-1$-dimensional hypersurfaces.
We proceed in an analogous way to the two-dimensional case above, again using the same pullback metric $g_t$.
A bounded $d-1$-dimensional hypersurface $S$ can be parameterised by a function $\gamma:[0,1]^{d-1}\to\mathbb{R}^d$.
To compute the $d-1$-dimensional area of the hypersurface $\Phi^t(S)$ one integrates $$\int_{[0,1]^{d-1}} \left(\det(\mathrm{Gram}_{g_t}(\mathbf{s}))\right)^{1/2}\ d\mathbf{s},$$ where 
$d\mathbf{s}=ds_1\,ds_2\cdots ds_{d-1}$ and $\mathrm{Gram}_{g_t}(\mathbf{s})_{ij}=\langle \partial\gamma(\mathbf{s})/\partial s_i,\partial\gamma(\mathbf{s})/\partial s_j\rangle_{g_t}$ is the $(d-1)\times(d-1)$ Gram matrix formed by partial derivatives of the parameterisation $\gamma$ with respect to the $d-1$ parameters.
In the case $d=2$, the Gram matrix is scalar-valued and we obtain \eqref{lineeqn2}.

\subsection{The dynamic Laplacian}
Consider a Riemannian manifold $(M,g)$ with volume form $\omega_g$.
The divergence of a vector field $X$ on $M$ is the unique scalar-valued function $\mathrm{div}(X):M\to\mathbb{R}$ satisfying $\mathrm{div}(X)\omega_g=X(\omega_g)$, where $X(\omega_g)$ is the Lie derivative of $\omega_g$ with respect to $X$.
The gradient $\nabla_g f$ is the unique vector field on $M$ satisfying $\langle \nabla_g f, X\rangle_g = X(f)$ for all vector fields $X$ on $M$, where $X(f)$ is the Lie derivative of $f$ with respect to $X$.
The Laplace--Beltrami operator with respect to a metric $g$ is defined as $\Delta_g f=\mathrm{div}(\nabla_{g}f)$.
The dynamic Laplacian \cite{F15} for $\Phi^t$ over the interval $[0,\tau]$ combines Laplace--Beltrami operators for the metrics $g_t$, $0\le t\le \tau$:
Let
\begin{equation}
\label{dynlap}
\Delta^D f:=\frac{1}{\tau}\int_0^\tau \Delta_{g_t}f\ dt,
\end{equation}
where $\Delta_{g_t}$ is the Laplace--Beltrami operator on $M$ with respect to the metric $g_t$.
Recall that each of the pullback metrics is a metric on $M$ and so \eqref{dynlap} is an average of Laplace--Beltrami operators, each of which is defined on $M$.

We now consider the (Neumann) eigenproblem for the dynamic Laplacian
\begin{equation}
\label{dynlapeig}
\Delta^D f=\lambda f,
\end{equation}
with $\int_0^\tau\langle \nabla_{g_t} f(x), n(x)\rangle_e\ dt=0$ for all $x\in\partial M$.
Because $\Delta^D$ is a linear combination of Laplace--Beltrami operators, it inherits the symmetry and ellipticity properties of these operators and therefore also many of their spectral properties.
The eigenspectrum of $\Delta^D$ is a countable sequence $0=\lambda_1^D\ge \lambda_2^D\ge\cdots$ of negative numbers diverging to $-\infty$, and the leading eigenfunction is constant.
The eigenfunctions $f_1^D, f_2^D,\ldots$ form a complete orthonormal set for $L^2(M)$ \cite{F15}.

\subsection{Dynamic Cheeger problems and dynamic Laplace eigenfunctions}
We now discuss the links between dynamic Cheeger problems and eigenfunctions of the dynamic Laplacian, in analogy to Section \ref{sec:cheelap}.
Recall the dynamic Cheeger ratio $H^D(A)$ for a set $A\subset M$ in equation \eqref{dyncheegerratio}.
When $H^D(A)$ is low, this means that the relative evolved boundary size of $\Phi^t(A)$ is low on average over the time interval $[0,\tau]$.
Given a dynamical evolution of a manifold $M$ under $\Phi^t, 0\le t\le \tau$, to find $n$ finite-time coherent sets we search for dynamic $n$-packings $\{A_1,\ldots,A_n\}$ of $M$ with low $H^D(\{A_1,\ldots,A_n\}):=\max_{1\le k\le n} H^D(A_k)$.
To quantify how well a choice can be made for a given finite-time dynamics we define the \emph{dynamic Cheeger constant} of $(\Phi^t,M)$ for $t\in[0,\tau]$ as\footnote{see \cite{F15} for $n=2$, and \cite{FR24} for general $n$.}
\begin{equation}
\label{dyncheegerconst}
h_n^D:=\inf\{H^D(\{A_1,\ldots,A_n\}): A_1,\ldots,A_n\subset M\mbox{ are disjoint}\}.
\end{equation}

We will use eigenfunctions of the dynamic Laplace operator to find good candidate dynamic $n$-packings.

\begin{theorem}[Theorem 3.19 \cite{FR24}]
\label{dynrockthm}
Consider an eigenpair $(\lambda^D,f^D)$ for the dynamic Laplacian $\Delta^D$.
\begin{enumerate}
\item If $f^D$ has at least $n$ nodal domains then $h^D_n\le\sqrt{-2\lambda^D}$.
\item For each $k=1,\ldots,n$, let $N_k\subset M$ denote the $k^{\rm th}$ nodal domain and let $A_k(c_k):=\{x\in N_k: f^D(x)^2\ge c_k\}\subset M$ denote the superlevel set of $(f^D)^2$ restricted to $N_k$.
    Then for a positive measure set of $(c_1,\ldots,c_n)\in\mathbb{R}^n$, one has 
    \begin{equation}
    \label{feassolD}
    H^D(A_1(c_1),\ldots,A_n(c_n))\le \sqrt{-2\lambda^D}.
    \end{equation}
\end{enumerate}
\end{theorem}
The interpretation of this theorem follows the discussion of Theorem \ref{ncheegerthm}.
By \eqref{dyncheegerratio}, part 1 says that the closer the eigenvalue $\lambda^D$ is to zero, the more easily we can decompose $M$ into $n$ disjoint sets $A_1,\ldots,A_n$ such that the average boundary size of each evolved family $\{\Phi^t(A_k)\}_{0\le t\le \tau}$ is small.
Part 2 provides an ansatz for identifying such families using the dynamic Laplacian eigenfunctions.

\begin{figure}
  \centering
  \includegraphics[width=\textwidth]{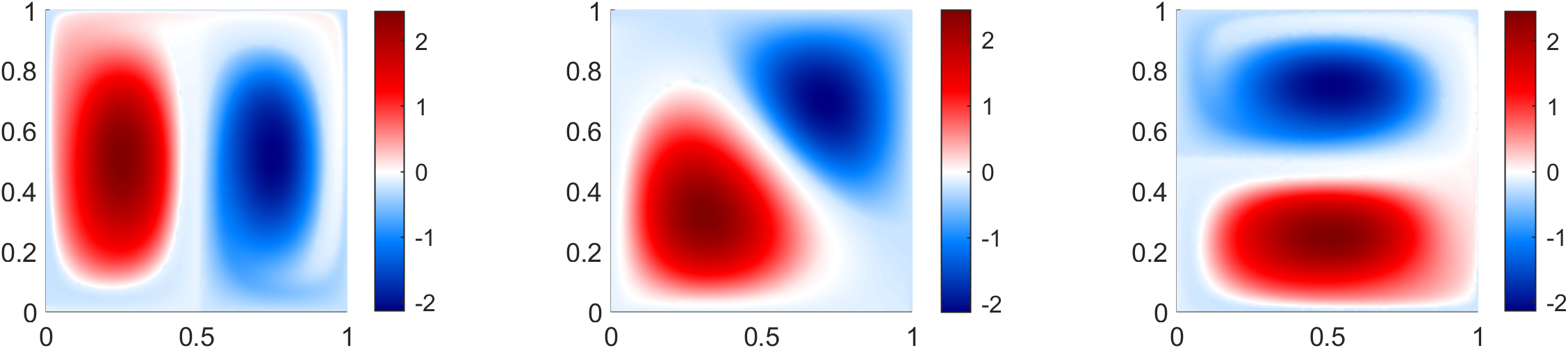}
  \caption{Left: Leading nontrivial eigenfunction $f_2^D$ for the rotating double gyre shown in Figure \ref{fig:meiss_triple}. Centre: Pushforward of $f_2^D$ by 1/2 a time unit to $t=1/2$, namely $f_2^D\circ(\Phi^{1/2})^{-1}$. Right: Pushforward of $f_2^D$ by 1 time unit to $t=1$, namely $f_2^D\circ(\Phi^{1})^{-1}$.}\label{fig:dgevecs}
\end{figure}

Figure \ref{fig:dgevecs}(left) shows the leading nontrivial eigenfunction $f_2^D$ of the dynamic Laplacian for the rotating double gyre dynamics discussed in Section \ref{sec:chaoticadvection}.
If we use superlevel sets of $(f^D_2)^2$ as per Part 2 of Theorem \ref{dynrockthm} we will obtain two sets $A_1, A_2$, which we can visualise as the deeper red and blue areas in Figure \ref{fig:dgevecs}(left).
We may push these red and blue regions forward under the dynamics $\Phi^t$ in the usual way, to obtain a function $f^D_2\circ(\Phi^t)^{-1}$ on $\Phi^t(M)$;  this is shown at $t=1/2$ and $t=1$ in Figure \ref{fig:dgevecs}(centre \& right).
These are snapshots of the evolving coherent families of sets $\Phi^t(A_1)$, $\Phi^t(A_2)$ for $t\in[0,\tau]$.
The red and blue structures that we see follow the dynamics visible in Figure \ref{fig:meiss_triple}, but the evolved boundaries of our sets $A_1, A_2$ are much smaller than the gray regions in Figure \ref{fig:meiss_triple}, and therefore much more coherent.

\subsection{Dirichlet boundary conditions}
\label{sec:dirichlet}
We have discussed the eigenproblem for the dynamic Laplace operator with Neumann boundary conditions.
These boundary conditions are particularly appropriate when the domains $\Phi^t(M)$ have a physically meaningful boundary and the finite-time coherent sets one seeks may touch the boundary. In situations where one explicitly wishes to exclude coherent sets that touch the boundary, or where the boundary of $\Phi^t(M)$ may be arbitrary because the manifolds $\Phi^t(M)$, $0\le t\le\tau$ live inside a larger ambient space, then Dirichlet boundary conditions -- i.e.\ the eigenfunctions are zero on $\partial M$ -- may be applied.
The main changes are \cite{FJ18,FR24}:
\begin{enumerate}
\item $\lambda^D_1$ is now strictly negative and $f^D_1$ may be chosen to be strictly positive, except on $\partial M$ where it vanishes.
\item Equation \eqref{dyncheegerratio} becomes
\begin{equation}
\label{dyncheegerratiodir}
h^D(A):=\frac{1}{\tau}\int_0^\tau \frac{|\partial(\Phi^t(A))|}{|\Phi^t(A)|}\ dt,
\end{equation}
so that the entire boundaries of $\Phi^t(A)$ contribute, even if part of these boundaries coincide with the boundaries of $\Phi^t(M)$.  With this adjustment Theorem~\ref{dynrockthm} holds for eigenpairs $(\lambda^D,f^D)$ satisfying Dirichlet boundary conditions.
\item Dirichlet boundary conditions force the superlevel sets $A_k(s_k)$ in Theorem~\ref{dynrockthm} to not touch the boundary of $M$.
\end{enumerate}

\section{Computation}

Our starting point will be several trajectories of a time-dependent dynamical system $x_{i,t_l}$, $0\le i\le N, 1\le l\le T$, where $x_{i,t_l}\in \Phi^{t_l}(M)$ is the location of the $i^{\rm th}$ trajectory at time $t_l$. 
These trajectories may be numerically generated from a known dynamical system $\Phi^t$ or they may be state observations made from some unknown dynamics.
The basic building block we need to approximate numerically is the Laplace--Beltrami operator for the pullback metric $g_t$.
Estimating $g_t=(\Phi^t)^*e$ numerically requires high-resolution data in order to accurately estimate the matrix $(D\Phi^t)^\top D\Phi^t$ (see Section \ref{sec:pullback}).
To avoid this high computational cost, we instead use the important relationship 
\begin{equation}
\label{transform}
\Delta_{g_t} f=\Delta_e (f\circ\Phi^{-t})\circ\Phi^t = (\Phi^t)^*(\Delta_e(\Phi^t_* f)).
\end{equation}

\begin{figure}
  \centering
  \includegraphics[width=\textwidth]{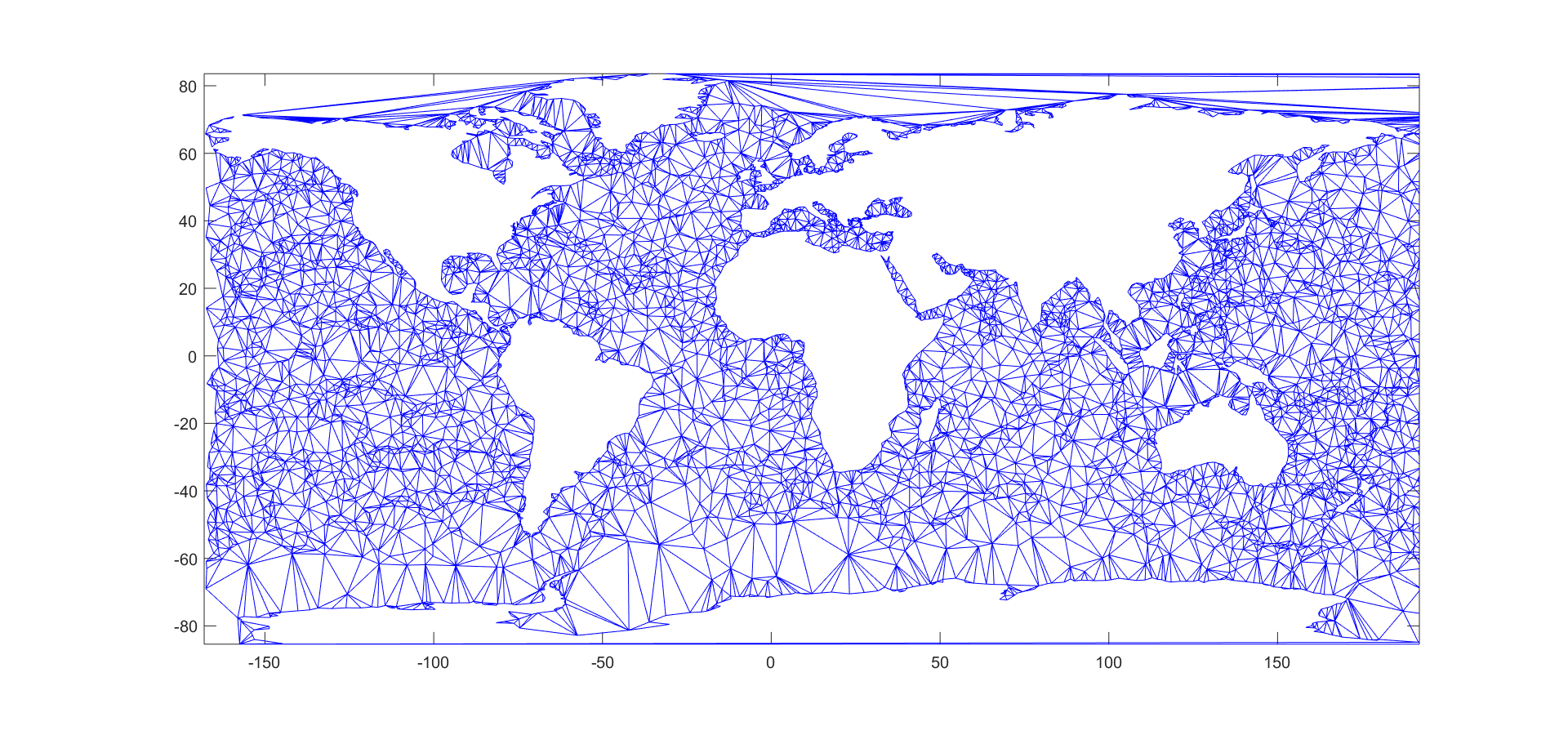}
  \caption{Delaunay mesh of approximately 3000 Argo float positions (nodes) at December 2013. Abernathey \emph{et al.} \cite{ABFS22} identified the major long-term coherent sets at 1000m depth by meshing Argo float locations $x_i^{t_l}$ at each month $t_l$ over six years.}\label{fig:oceanmesh}
\end{figure}

Because we wish to estimate Laplace--Beltrami operators, a natural established numerical framework is the finite-element method.
This has been adapted to the dynamic Laplacian in \cite{FJ18}.
We wish to solve the eigenproblem $\Delta^D f=\lambda f$ with the dynamic Neumann boundary condition.
We begin by integrating both sides of this equality against a test function $\tilde{f}$ using the volume measure $V_e$ for the metric $e$, i.e.\ Lebesgue measure.
This yields the weak form of the dynamic Laplacian eigenproblem, namely 
\begin{equation}
\label{weakeig}\int_M (\Delta^D f)\tilde{f}\ dV_e=\lambda\int_M f\tilde{f}\ dV_e\mbox{ for all test functions $\tilde{f}$.}
\end{equation}
The weak form has the immediate advantage of eliminating the Neumann\footnote{Dirichlet boundary conditions are also eliminated by simply restricting the functions $f$ to those that vanish at the boundary of $M$.} boundary condition \cite{FJ18}.
By \eqref{dynlap} we see that the LHS of \eqref{weakeig} is a time integral of $\int_M (\Delta_{g_t}f)\tilde{f}\ dV_e$ over the time interval $[0,\tau]$. 
Taking one of these terms and applying \eqref{transform} we obtain
$$\int_M (\Phi^t)^*(\Delta_e(\Phi^t_* f))\tilde{f}\ dV_e.$$
Changing variables with $\Phi^t(M)$ and performing integration by parts yields
\begin{equation}
\label{futureweak}-\int_{\Phi^t(M)} \nabla(\Phi^t_* f)\cdot \nabla(\Phi^t_*\tilde{f})\ dV_e.
\end{equation}
We now introduce the FEM discretisation, which is a basis of nodal hat functions $\{\varphi^{t_l}_i\}_{j=1}^N$ defined on a mesh of $\Phi^t(M)$ with $N$ nodes $\{\Phi^{t_l}x_i\}_{i=1}^N$, namely the trajectories we were initially supplied with.
As a real-world example, Figure \ref{fig:oceanmesh} shows a mesh of around 3000 deep-ocean ARGO floats at December 2013.
Each $\varphi^{t_l}_i$ is locally supported on mesh elements containing the node $\Phi^{t_l}x_i$ as a vertex, and $\varphi^{t_l}_i(\Phi^{t_l}x_i)=1$, $i=1,\ldots,N$.
We form a (sparse) matrix $D^{t_l}_{ij}$ by evaluating (\ref{futureweak}) with $\Phi^t_*f=\varphi^{t_l}_i$ and $\Phi^t_*\tilde{f}=\varphi^{t_l}_j$, $1\le i,j,\le N$:
\begin{equation}
\label{futureweak2}
D^{t_l}_{ij}:=-\int_{\Phi^{t_l}(M)} \nabla\varphi^{t_l}_i\cdot \nabla\varphi^{t_l}_j\ dV_e,
\end{equation}
which is a standard FEM integral.
We have thus handled the left hand side of \eqref{weakeig}, obtaining the matrices $D^{t_l}_{ij}$ for $t$ in a discrete set of points $\{t_1,\ldots,t_T\}$ sampled from the time interval $[0,\tau]$. 

For the right hand side of \eqref{weakeig} we  substitute $f=\varphi^0_i, \tilde{f}=\varphi^0_j$ to obtain the $N\times N$ sparse matrix $M_{ij}:=\int_M \nabla\varphi^0_i\cdot\nabla\varphi^0_j\ dV_e$.
One now solves the discretised problem 
\begin{equation}
\label{femeqn}
\left(\frac{1}{T}\sum_{l=1}^T D^{t_l}\right)v=\lambda Mv
\end{equation} to find eigenvalues $\lambda$ and corresponding basis coefficient vectors $v$, and reconstructs eigenfunctions on $M$ as $f=\sum_{j=1}^N v_j\varphi^0_j$.

This FEM approach is very robust and accurate, and is numerically efficient for dimensions 2 and 3; it has been used for all eigenfunction images in this article.
It has been tested in many situations on sparse trajectory data \cite{FJ18,FRS19,ABFS22} and is insensitive to highly varying densities of trajectory data sampling (see e.g.\ Figure 10 \cite{FJ18}).
The method also copes very well with situations where observed trajectories are incomplete and many observations are missing.
Figure 12 \cite{FJ18} illustrates this for the rotating double gyre eigenfunction in Figure \ref{fig:dgevecs}.

\begin{figure}
  \hspace*{-1.2cm}\includegraphics[width=1.2\textwidth]{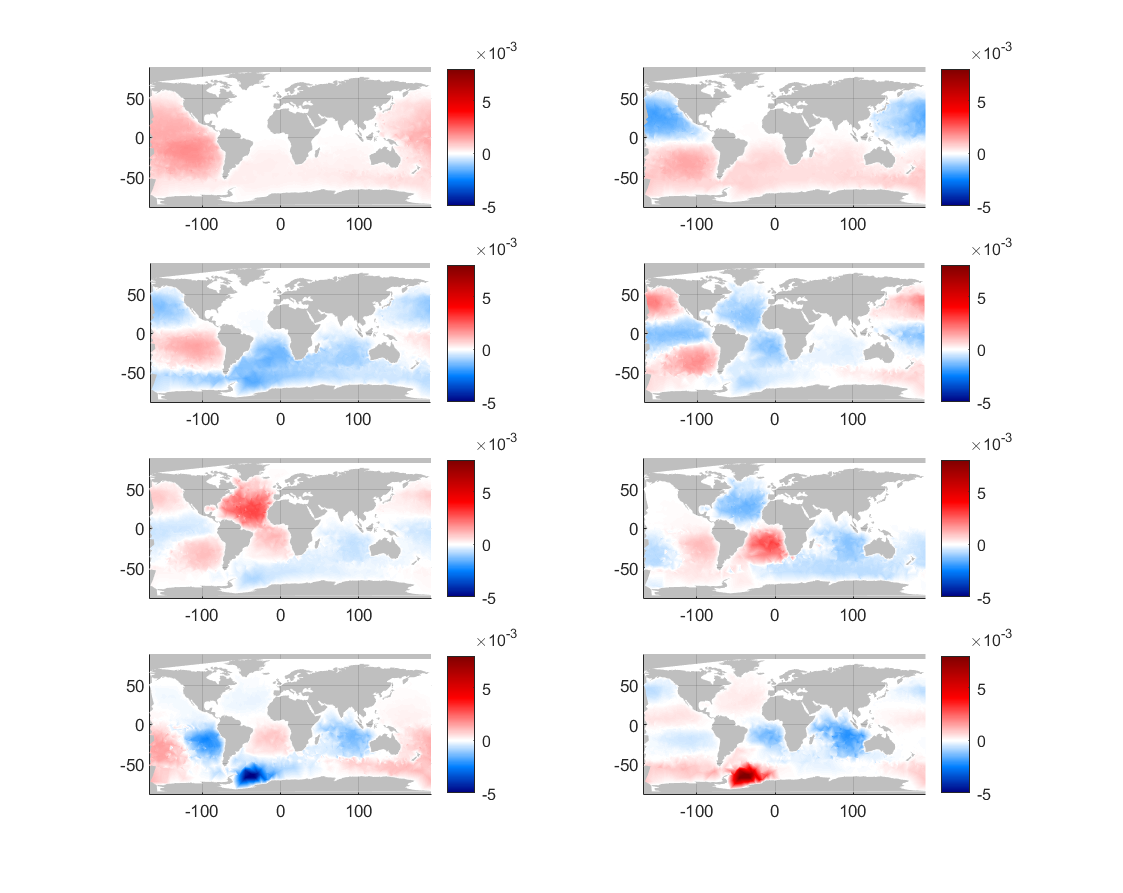}
  \caption{Leading eight eigenfunctions $f_1^D,\ldots,f_8^D$ of the dynamic Laplacian computed from Argo trajectory
data over the six years 2011-2017. The eigenfunctions are represented as scalar fields at January 2011. Image from \cite{ABFS22}.}\label{fig:ocean1}
\end{figure}

The identification of coherent regions in the global ocean at 1000m depth from Argo drifter observations is an example of a  real-world computation requiring a method than can handle sparsity, varying data density, and missing data.
Figure \ref{fig:oceanmesh} illustrates the sparsity and the varying data density.
Figure \ref{fig:ocean1} displays eigenfunctions of the dynamic Laplacian built from around 8000 individual deep-sea Argo float trajectories, where only 3000 floats are providing location signals in any given month, and only 10\% of floats continue to transmit over the six years from 2011--2017.
In Figure \ref{fig:ocean1} the dynamic Laplacian was estimated in a single computation constructed from the many partial and non-overlapping-in-time trajectories of Argo floats.

Finally, we note that there are no free parameters to choose and the method preserves the symmetric structure of the original eigenproblem.
MATLAB code and demos are available at \url{https://github.com/gaioguy/FEMDL}.
Alternate FEM-based strategies and generalisations to weighted manifolds are discussed in detail in \cite{FJ18}.

\subsection{Alternate approaches to estimate $\Delta^D$ and related objects}

\paragraph{Ulam and finite-difference:}
In \cite{F15} the pushforwards and pullbacks $\Phi_*^t$ and $(\Phi^t)^*$ in \eqref{transform} were estimated with Ulam-discretisations and the Laplacian $\Delta_e$ was discretised by finite-difference.
This approach was generalised to the measure-preserving setting on weighted manifolds in \cite{FK20}, and we refer the reader to \cite{F15,FK20} for a full explanation of this approach. A large number of points were used in \cite{F15,FK20} to estimate the Ulam matrices representing $\Phi_*^t$ and $(\Phi^t)^*$, but similar results can be obtained with many fewer points. Because this approach is grid-based, a reasonably uniform sampling of initial conditions is needed.

\paragraph{Radial basis function collocation:}
The paper \cite{FJ15} replaces the indicator function bases of the Ulam approach with smooth radial basis functions and proceeds to approximate $\Phi_*^t$, $(\Phi^t)^*$ and $\Delta^D$ by collocation. The shape parameter -- which determines the width of the basis functions -- needs to be selected and the method requires a suitable positioning of basis function centres.

\paragraph{Related objects:}
We noted earlier that the dynamic Laplacian is a type of spectral clustering on manifolds.
 A lower-order approach, compared to FEM, is to construct a graph Laplacian based on initial conditions of trajectories \cite{hadjighasem}, where edges are weighted according to a dynamic distance introduced in \cite{FPG15}. 
In a similar direction, \cite{padbergschneide} apply constant edge weights that record if one trajectory comes within a threshold of another. 
In these approaches, cutoffs for weights or distances need to be selected.
Varying point density can influence the results, whereas this is automatically adjusted for in the FEM approach, which can itself be interpreted as a graph Laplacian, see the discussion in Section 3.4 \cite{FJ18}.

\section{Finding multiple coherent sets: SEBA}
\begin{figure}
  \centering
  \includegraphics[width=\textwidth]{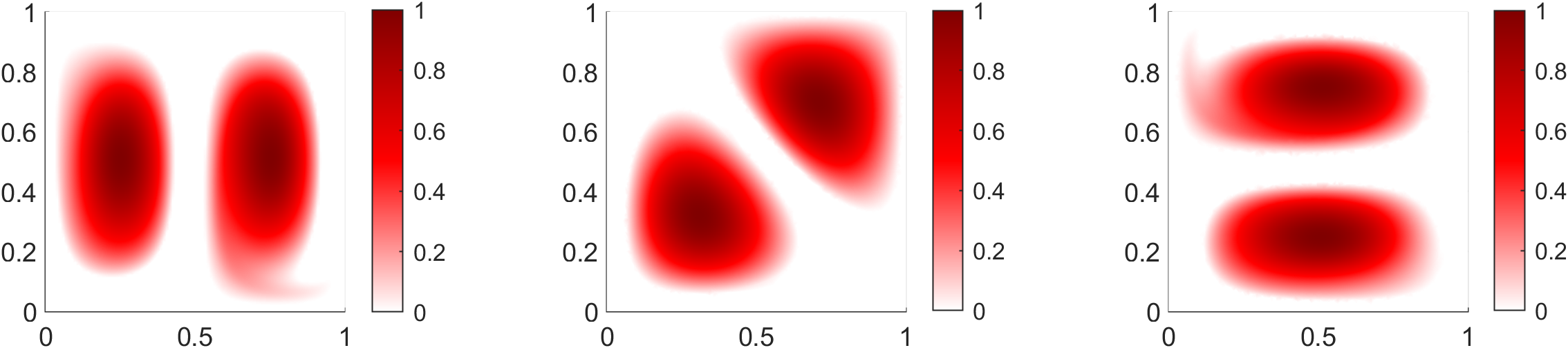}
  \caption{Left: Sum of two SEBA vectors arising from SEBA applied to the dynamic Laplacian eigenfunctions $f^D_1\equiv\mathbf{1}$ and $f^D_2$ from Figure \ref{fig:dgevecs}(left) for the rotating double gyre. Centre: Pushforward of the left panel by 1/2 a time unit with $\Phi^{1/2}_*$ to $t=1/2$. Right: Pushforward of the left panel by 1 time unit with $\Phi^{1}_*$ to $t=1$.}\label{fig:dgseba}
\end{figure}

In Section \ref{sec:seba1} we discussed how the SEBA algorithm could create sparse linear combinations of Laplace eigenfunctions to automatically separate individual nodal domains and create additional candidate sets for $n$-packings.
The same methodology can be applied to the eigenfunctions $f^D_1,\ldots,f^D_K$ of the dynamic Laplace operator to create a sparse basis $\{s^D_1,\ldots,s^D_K\}$ whose span approximates $\spn\{f^D_1,\ldots,f^D_K\}$.
Each sparse function $s^D_k$ supports a single coherent set and in this way we disentangle individual coherent sets from multiple eigenfunctions.
Performance guarantees for candidate coherent sets $A_k(c_k):=\{x\in M: s^D_k(x)\ge c_k\}$ are available in direct analogue to Part 2 of Theorem \ref{dynrockthm};  see \cite{FR24}.
By applying SEBA to the the leading eigenfunctions $f^D_1\equiv\mathbf{1}$ and $f^D_2$ from Figure \ref{fig:dgevecs}(left), we obtain two (sparse) SEBA vectors $s^D_1, s^D_2$ taking values in $[0,1]$.
The sum of these two vectors is shown in red in Figure \ref{fig:dgseba}(left) and their pushforwards under $\Phi^t$ shown in the centre and right panels.

\begin{figure}
  \hspace*{-1.2cm}\includegraphics[width=1.2\textwidth]{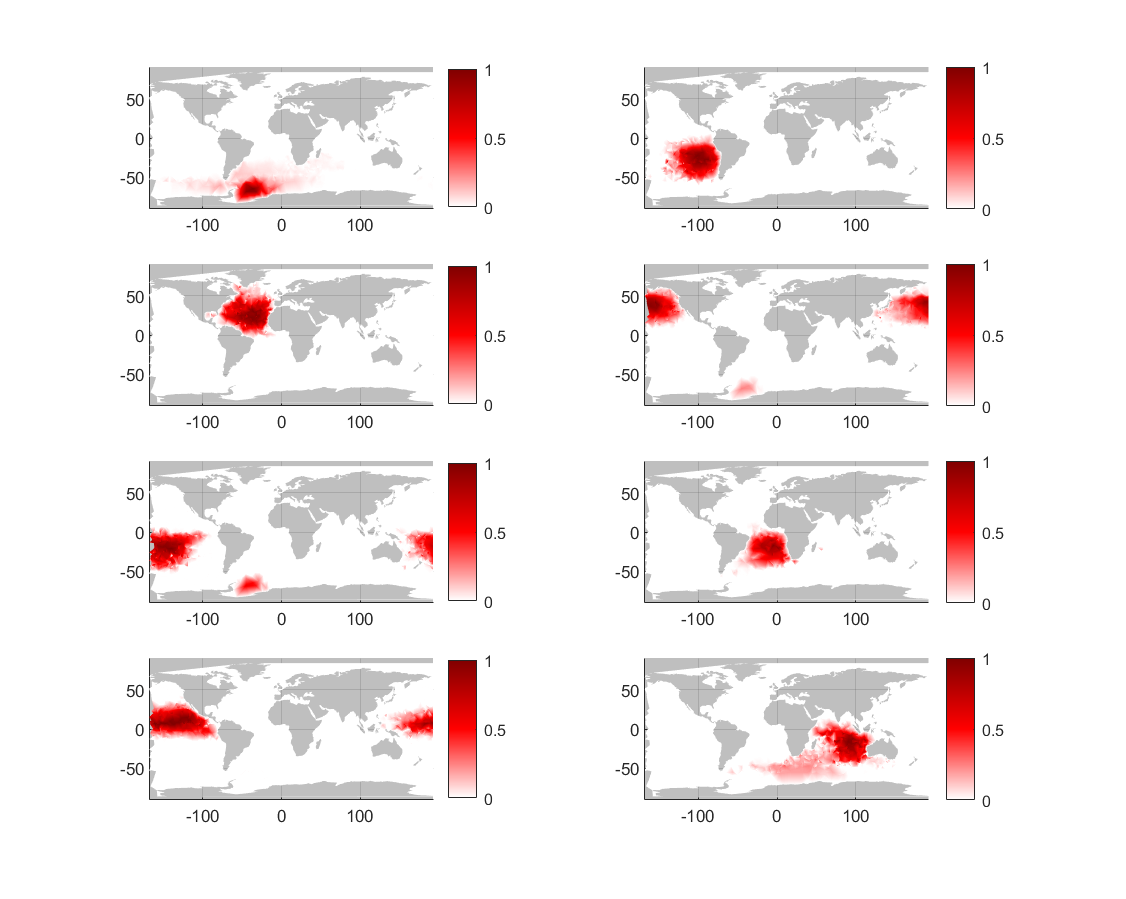}
  \caption{Result of SEBA applied to the eight eigenfunctions shown in Figure \ref{fig:ocean1} to construct eight SEBA vectors $s_1^D,\ldots,s_8^D$. Each SEBA vector is supported on an individual coherent region of the ocean at January 2011, each of which remains coherent under ocean currents at 1000m depth over a six-year duration from 2011-2017. Image from \cite{ABFS22}.}\label{fig:sebaocean}
\end{figure}

Similarly, we can apply SEBA to the eight eigenfunctions shown in Figure \ref{fig:ocean1} to separate the eight dominant coherent regions in the global ocean at 1000m depth;  see Figure \ref{fig:sebaocean}.
With the strongest coherent regions clearly separated, it has been noted in \cite{ABFS22} that they are largely aligned with the major ocean boundaries, but with a bias toward the
east of the basin due to the stronger-mixing currents on their western boundaries.
The Weddell gyre off the Antarctic coast is  highlighted.
The North and South
Atlantic is separated by the equatorial countercurrent. The subpolar and subtropical gyres in the North Pacific are separated 
and the South Pacific Gyre splits into eastern and western
components.
The reader is referred to \cite{ABFS22} for more detailed analysis. 
These results are broadly consistent with coarser results \cite{FPG15,BK16} using floating drifters on the ocean surface.

Sometimes the trajectory data is such that it does not support exactly $K$ coherent sets in the leading $K$ eigenfunctions.
We have already seen this issue in Section \ref{sec:seba1} where the eigenfunctions in Figure \ref{fig:blob} supported decompositions into $2, 4,$ and $5$ regions, but not $3$ regions.
The SEBA vectors themselves contain a simple ``reliability'' test, namely the magnitude of the most negative value;  see Sections 4.1 and 4.2 \cite{FRS19}). 
A discussion of assessing spectral gaps is in Section 4.2 \cite{FRS19}.

\section{Weights, dynamic metrics, and zero-noise limits} 

\paragraph{Measure-preserving dynamics on weighted Riemannian manifolds.} We have mainly discussed the special situation where we have a flat manifold $M$ embedded in Euclidean space and $\Phi^t$ preserves Euclidean volume.
We may also consider a dynamic Laplacian on a sequence of weighted Riemannian manifolds $(M_t,h_t,m_t)$, where $M_t=\Phi^t(M)$, $m_t$ is a Riemannian metric on $M_t$ and $h_t:M_t\to\mathbb{R}^+$ is a positive weight function.
The measure of $A\subset M$ is $\int_{A} h_0\ dV_{m_0}$, where $V_{m_0}$ is the volume measure induced by the metric $m_0$.
The dynamics $\Phi_t$ can be made \emph{measure preserving} by defining positive weight functions $h_t$ so that the measure of a set $A\subset M$ equals the measure of $\Phi^t(A)$ on $M_t$: $\int_{\Phi^tA} h_t\ dV_{m_t}=\int_A h_0\ dV_{m_0}$, for $0<t\le \tau$.
One now defines 
\begin{equation}
\label{dynlapweighted}
\Delta^Df:=\frac{1}{\tau}\int_0^\tau \Delta_{h_t,g_t}f\ dt,
\end{equation}
where $\Delta_{h,g}f=\frac{1}{h}\mathrm{div}(h\nabla_g f)$ is the $h$-weighted Laplace--Beltrami operator for $g$.

In the context of fluid flow, this weighted setting is useful for finding finite-time coherent sets in compressible fluid flows, where Euclidean volume is not preserved, but fluid mass is conserved.
The weight functions $h_t$ track the evolution of an initial fluid density $h_0$ (with respect to $V_{m_0}$) under $\Phi^t$.
This weighted setup also covers the situation where one is interested in coherence of a passive tracer carried by the flow and this passive tracer is initialised with a spatial concentration different to the fluid density; e.g.\ in the volume-preserving case, the fluid density would be uniform.
One sets $h_0$ to be the initial density of the passive tracer, computed with respect to $V_{m_0}$, and future passive tracer densities are tracked by the weight functions $h_t$.
We refer the reader to \cite{FK20} for the relevant constructions and \cite{FJ18} for numerical computations.

 
\paragraph{Is the dynamic Laplacian a Laplace--Beltrami operator?}

One may ask whether there is a dynamic Riemannian metric $g^D$ so that $\Delta^D=\Delta_{g^D}$.
Let us return to the simplified setting of a flat manifold $M\subset\mathbb{R}^d$ endowed with the Euclidean metric $e$ and $\Phi^t$ Lebesgue volume preserving; i.e.\ $V_{g_t}=V_e$ for all $t\in[0,\tau]$.
The appropriate construction for $g^D$ is to set $g^D=\left(\frac{1}{\tau}\int_0^\tau g_t^{-1}\ dt\right)^{-1}$, where $g^{-1}$ represents the dual metric, which is a symmetric, positive-definite, bilinear form acting on pairs of covectors in cotangent space \cite{KK,KS}.
One then obtains that $\Delta^D$ is the weighted Laplace--Beltrami operator $\Delta_{\theta, g^D}$, where $\theta$ is the Radon-Nikodym derivative $dV_e/dV_{g^D}$. 
This operator acts on the weighted manifold $(M,\theta,g^D)$, where the role of the weight $\theta$ is to transform the weighted volume of a set $A\subset M$, namely $\int_A \theta\ dV_{g^D}$, into the Lebesgue volume of A, $\int_A \mathbf{1}\ dV_{e}$.
Similar results hold for the measure-preserving case in the previous paragraph.  See \cite{KK,KS,FL24} for details.

\paragraph{The dynamic Laplacian as a zero-noise limit.}

In Section \ref{sec:chaoticadvection} we discussed the connection between dynamic isoperimetric problems and finite-time coherent sets by considering advective dynamics with small-noise diffusive dynamics.
In fact, there is a formal connection in the following sense. 
Finite-time coherent sets were initially introduced and motivated by minimally mixing behaviour of subsets in the presence of small noise \cite{FSM10,F13}.
The classical transfer operator approach \cite{F13} uses leading singular vectors of a diffusive transfer operator $\mathcal{L}_\epsilon = \mathcal{D}'_\epsilon\circ \Phi^\tau_*\circ\mathcal{D}_\epsilon$, where $\mathcal{D}'_\epsilon$ and $\mathcal{D}_\epsilon$ are isotropic diffusion (of scale $\epsilon$) operators on $M$ and $\Phi^\tau(M)$, respectively.
One may ask what happens in the $\epsilon\to 0$ limit.
The naive approach of considering leading eigenvectors of $\lim_{\epsilon\to 0} \mathcal{L}_\epsilon^*\mathcal{L}_\epsilon$ is uninformative because this limit is the identity operator.
However, a finer analysis \cite{F15} shows that the dynamic Laplacian emerges in an appropriate limit.
For $f\in C^3(M)$
$$\lim_{\epsilon\to 0}\frac{(\mathcal{L}_\epsilon^*\mathcal{L}_\epsilon-I)f}{\epsilon^2}=C\cdot(\Delta_e + (\Phi^\tau)^*(\Delta_e(\Phi_*^tf))=C\cdot\left(\Delta_{g_0} + \Delta_{g_\tau}\right)/2,$$
for a constant $C$ depending on the common $\epsilon$-scaled kernel of the diffusion operators $\mathcal{D}_\epsilon$.
The RHS above is a dynamic Laplace operator where the average of Laplace--Beltrami operators is taken at the end points of the time interval $[0,\tau]$.
This result was proved for general weighted Riemannian manifolds \cite{FK20} and the limit shown to be uniform in the $C^0$-topology for $f$ ranging over $C^3(M)$.
Related calculations are in \cite{BK16} and \cite{KS}.

\section{When finite-time coherent sets come and go:  the inflated dynamic Laplacian}

The dynamic Laplacian tracks dynamics over a time interval $[0,\tau]$ and is an average of Laplace--Beltrami operators $\Delta_{g_t}$ with respect to metrics that are pullbacks under the dynamics.
Because of this averaging property, the dynamic Laplacian succeeds in constructing dynamic isoperimetric properties of sets, and consequently finite-time coherent sets for the interval $[0,\tau]$.
In several real-world situations, one may have complex dynamics in which many finite-time coherent sets emerge and disappear, each with their own distinct lifetimes. 
The dynamic Laplacian may be deployed on distinct time intervals to either separately find these sets and/or stitch distinct computations together across time \cite{schneide2022}.
A recent alternative is to create an
\emph{inflated dynamic Laplace operator} on the spacetime manifold $[0,\tau]\times M$ for a perhaps larger $\tau$ than many of the lifetimes of the coherent sets.
On each time fibre $\{t\}\times M$, one applies $\Delta_{g_t}$ (acting in space), while the distinct time fibres are tied together with a one-dimensional Laplace operator on $[0,\tau]$.
The eigenfunctions of this inflated dynamic Laplacian encode the appearance and disappearance of multiple finite-time coherent sets, whose individual lifetimes can be extracted by applying the SEBA algorithm in spacetime.
Much of the theory of the dynamic Laplacian can be naturally modified for the inflated dynamic Laplacian;  see \cite{FK23,AFK24} for details.

\subsection{Clusters in time-varying networks}

The dynamic Laplacian and the inflated dynamic Laplacian have counterparts on time-varying networks.
For example, one may consider permutation dynamics on the nodes of the network and develop a corresponding dynamic spectral theory for dynamic graph cuts \cite{FK15}.
The concept of a supra-Laplacian \cite{Gmez2013,sole2013spectral} is related to a spacetime network version of the inflated dynamic Laplacian and recent work \cite{FKK24} has shown how spacetime spectral clustering in time-varying networks can automatically detect the appearance and disappearance of clusters.

\section{Acknowledgements}
The author thanks Yat-Long Lee for a careful reading of the manuscript. GF's research is partially supported by an Australian Research Council Discovery Project.  

\bibliographystyle{spmpsci_unsrt}

 \end{document}